\documentclass[]{amsart}
 \usepackage[]{amsmath, amsthm, amsfonts}
 \usepackage[all]{xy}
 \usepackage{hyperref}
 \usepackage{amssymb}

 \newcommand {\C} {{\mathbb C}}

 \newcommand {\Q} {{\mathbb Q}}
 \newcommand {\PP} {{\mathbb P}}

 \newcommand {\dt} {{\bullet}}

 \newcommand {\mU}{\mathcal{U}}

 \newcommand {\mR}{\mathcal{R}}
 \newcommand {\mW}{\mathcal{W}}
 \newcommand {\T}{\mathcal{T}}
 
 \newcommand {\tS}{\tilde{S}}

 \newtheorem{thm}[subsection]{Theorem}
 \newtheorem{cor}[subsection]{Corollary}
 \newtheorem{lemma}[subsection]{Lemma}

\begin{document}

 \title{Functoriality of the coniveau filtration}

 \author{Donu Arapura}\thanks{First author partially supported by the
   NSF} \author{ Su-Jeong Kang}
 \address{Department of Mathematics\\
   Purdue University\\
   West Lafayette, IN 47907\\
   U.S.A.}  \maketitle

 \begin{abstract}
  It is shown that the coniveau filtration on the cohomology 
of smooth projective varieties is preserved up to shift
 by  pushforwards, pullbacks and  products.
 \end{abstract}

 There is a natural descending filtration on the singular cohomology
 of a complex smooth projective variety called the coniveau
 filtration.  A class lies in the $p$th level of coniveau if it is
 supported on a closed subvariety of codimension at least $p$, that is
 if it vanishes on the complement of such a subvariety. The generalized
 Hodge conjecture  would imply, rather trivially, that the
 coniveau filtration is compatible with pushforwards, pullbacks and
 products. The purpose of this paper is to prove this statement
 unconditionally. 

 The main difficulty is proving the compatibility of coniveau with
 respect to pullbacks along closed immersions. Given a closed
 immersion $i:X\to Y$ and a class $\alpha\in H^k(Y)$ supported on
 subvariety $S$, $i^*\alpha$ is supported on $S\cap X$. So the only
 issue is what to do when  $S$ and $X$ do not intersect properly. The
 classical solution to this kind of problem is to prove a
 moving lemma. However, our situation is fairly rigid which makes
 this approach difficult. We solve the problem by blowing up $Y$ and
 choosing suitable lifts of $X$ and $S$ in general position.  Much of
 the work in this paper involves showing that this suffices.

 \section{Main theorem}

We use the symbol $\subset$ for nonstrict inclusion.
 All our varieties will be defined over $\C$. Given a variety $X$, we
 write $H^i(X)$ for its singular cohomology with rational
 coefficients.  We will say that a class $\alpha\in H^i(X)$
is supported on a closed subvariety $S$ if it lies in the
kernel 
$$\text{Ker}[H^i(X)\to H^i(X-S)]$$
By \cite[Proposition 8.2.8]{Deligne3}, this is equivalent to
$\alpha$ lying in the image of a Gysin homomorphism. More precisely,
that
$$\alpha\in \text{Im}[H^{i-2q}(\tilde S)\to H^i(X)]$$
where $q= \text{codim}(S,X)$ and $\tS\to S$ is a desingularization.
The {\em coniveau} filtration on $H^i(X)$ is given by
 $$
 N^p H^i(X) = \sum_{codim S\ge p} \text{Ker}[H^i(X)\to H^i(X-S)]$$
 We call an element 
$\alpha\in N^pH^i(X)$ {\em irreducible} if it is supported on
one irreducible subvariety $S$.

The  generalized Hodge conjecture (GHC) predicts that $N^pH^i(X)$ is
the  maximal Hodge substructure of $F^pH^i(X)$ \cite{groth2,Lewis}.
The  filtration on $H^i(X)$ by maximal Hodge structures in 
$F^{\bullet}H^i(X)$ is easily seen to have good
functorial properties. Our main result is that this so for coniveau without
assuming GHC.

\begin{thm}\label{thm:functorialityofN}
  The coniveau filtration $N^\dt$  is preserved (up to shift)
by pushforwards,  exterior products and pullbacks. More precisely;
  \begin{enumerate}
  \item If $f:X\to Y$ is a map of smooth projective varieties of
    dimensions $n$ and $m$ respectively, then
    $$f_*(N^pH^i(X))\subset N^{p+m-n}(H^{i+2(m-n)}(Y))$$
  \item $$N^p(H^i(X))\otimes N^q(H^j(Y))\subset
    N^{p+q}H^{i+j}(X\times Y)$$

  \item If $f$ is as above, then $$f^*(N^pH^i(Y))\subset N^pH^i(X)$$
  \end{enumerate}
\end{thm}

This was stated in \cite{Arapura}, but the proof there was incomplete.

\begin{cor}
    $$
    N^p H^i(X) \cup N^q H^j(X) \subset N^{p+q}H^{i+j}(X) $$
  \end{cor}

  \begin{proof}
    The cup product is a composition of exterior product and
restriction to the diagonal.
  \end{proof}

We give the proof of the first two parts of the theorem now,
since we will need them.
The proof of the last part will be postponed until the final
section.

\begin{proof}[Proof of parts 1 and 2]
  An irreducible element $\alpha\in N^pH^i(X)$ lies in the image of a map
  $k_*(H^{i-2q}(T))$ where $k:T\to X$ is a morphism from a smooth
  projective variety of dimension $n-q\le n-p$.  Therefore
  $$f_*(\alpha)\in (f\circ k)_*H^{i-2q}(T) \subset
  N^{p+m-n}(H^{i+2(m-n)}(Y))$$
  This proves the first part.

  For the second statement. Let $\alpha\in N^pH^i(X)$ and 
$\beta\in N^qH^j(X)$ be irreducible classes. Then $\alpha$ and
$\beta $ are supported on subvarieties 
 $S\subset X$ and $T\subset Y$  such that $\dim S\le \dim X-p$ and
  $\dim T\le \dim Y-q$. Then $\dim (S\times T) \le \dim (X\times Y)
  -p-q$.  It follows that $\alpha\times \beta$, which is supported on
 $S\times T$, lies in  $N^{p+q}(H^{i+j}(X\times Y))$ as expected.

\end{proof}

\section{Lemmas}

In this section, we give the  basic lemmas needed to finish the proof
of the theorem. Many of these lemmas are nothing but special cases of it.

\begin{lemma}\label{lemma:surjective}
  Let $f:X \to Y$ be a surjective map of smooth projective varieties,
  then $f^*(N^p H^i(Y)) \subset N^pH^i(X)$.
\end{lemma}
\begin{proof}
Suppose that $\alpha\in N^pH^i(Y)$ is an irreducible class supported
on an irreducible variety
   $S\subset Y$ of codimension $q\ge p$.
 The preimage $f^{-1}S$ will have codimension less than
  or equal to $q$. By taking general hyperplane sections, we can find
  a cycle $Z\subset f^{-1}S$ of codimension exactly $q$ surjecting
  onto $S$. By stratification theory \cite[pp. 33-43]{gm}, we can find
  a proper Zariski closed set $Z''\subset Z$ containing the union of
  singular loci $Z_{sing}\cup f^{-1}S_{sing}$, such that the map
  $f:X-Z''\to Y-f(Z'')$ is locally trivial along tubular neighborhoods
  of $Z'=Z-Z''$ and $S'=S-f(Z'')$.  To make the last condition
  precise, consider the diagram $$
  \xymatrix{
    N_{Z'}\ar[rr]\ar[rd]^{\pi}\ar[dd]_{f} &  & X'\ar[dd]^{f} \\
    & Z'\ar[ru]^{k}\ar[dd]^{f}\ar@/^/[lu] &  \\
    N_{S'}\ar[rr]\ar[rd]^{\pi} &  & Y' \\
    & S'\ar[ru]^{i}\ar@/^/[lu] & } $$
  where $X' = X - Z''$, $Y'= Y -
  f(Z'')$, and $N_{Z'}$ $N_{S'}$ denotes appropriately chosen tubular
  neighbourhoods of $Z'$ and $S'$ respectively.  The above condition
  is that $N_{Z'}\to f^*N_{S'}$ is a locally trivial map of locally
  trivial fibre bundles (for the classical topology) over $Z'$.
  Fibrewise, we have an open immersion of $2q$ real dimensional
  oriented manifolds, and this induces an isomorphism of compactly
  supported $2q$ dimensional cohomologies.  Thus the Thom class
  $\tau_{Z'}$ of $N_{Z'}$, which can be viewed element of relative
  cohomology $H^{2q}_{Z'}(N_{Z'} ) = H^{2q}(N_{Z'}, N_{Z'}-Z')$,
  coincides with the pullback of the Thom class $\tau_S'$ on
  $f^*N_{S'}$.  The Gysin map $H^{i-2q}(Z')\to H^{i}(X')$ is given by
  $\alpha\mapsto \pi^*\alpha\cup\tau_{Z'}$ extended by $0$ to $X'$. A
  similar description holds for $(S',Y')$. It follows that we have a
  commutative diagram  $$
  \xymatrix{
    H^{i-2q}(S')\ar[r]\ar[d] & H^{i}(Y')\ar[d] \\
    H^{i-2q}(Z')\ar[r] & H^{i}(X') } $$
  Therefore, since  $$\alpha \in
  \text{Ker}[H^i(Y)\to H^i(Y-S)= H^i(Y'-S')],$$
  then its image in
  $H^i(X)$ maps to $$\text{Ker}[H^i(X)\to H^i(X-Z)=H^i(X'-Z')].$$
  This
  implies that $f^*$ preserves $N^p$.
\end{proof}

In the special case where $f$ is also flat there is a much simpler
proof, which we feel compelled to give.

\begin{proof}[Proof for flat maps]
  Let $\alpha\in N^pH^i(Y)$ be irreducible. 
Then it supported on   an irreducible subvariety $S$ of $Y$ with
  $\text{codim}(S,Y)=q\geq p$. Let $T = X
  \times_Y S$. Then $\text{codim}(T, X)=q \geq p$ and we have
  a commutative diagram $$
  \xymatrix{
    H^i(Y) \ar[r]^{\iota^* \quad} \ar[d]_{\pi^*} & H^i(Y-S) \ar[d]^{\pi^*} \\
    H^i(X) \ar[r]_{\iota^*_1 \quad} & H^i(X -T)
  } $$
  where $\iota:Y-S \hookrightarrow Y$, $\iota_1 : X - T
  \hookrightarrow X$ inclusions. Since $\alpha \in
  \text{Ker}(\iota^*)$, we have
  $$
  0=\pi^* \circ \iota^*(\alpha)=\iota^*_1 \circ \pi^*(\alpha) $$
  Hence, $$
  \pi^*(\alpha) \in \text{Ker}[H^i(X) \to H^i(X - T)]
  \subset N^pH^i(X) $$
\end{proof}

\begin{lemma}\label{lemma:genericallyfinite}
  Let $f:X \to Y$ be a generically finite map of smooth projective
  varieties, then $f^* : H^i(Y) \to H^i(X)$ is injective. Furthermore,
  $$
  f^* \left(H^i(Y) \right) \cap N^pH^i(X) = f^* \left( N^pH^i(Y)
  \right) $$
\end{lemma}
\begin{proof}
  Let $d=\text{deg}~f$. Then the first statement follows from $f_*
  f^*=d \cdot \text{id}$, where $\text{id}:H^i(Y) \to H^i(Y)$ the
  identity map. For the second statement, first note that $f^*$
  preserves coniveau by Lemma~\ref{lemma:surjective}, so it is enough
  to show $$
  f^* \left(H^i(Y) \right) \cap N^pH^i(X) \subset f^*
  \left( N^pH^i(Y) \right) $$
 Suppose that  $\alpha \in H^i(Y)$ satisfies
 $$
  f^*(\alpha) \in
  f^*\left(H^i(Y)\right) \cap N^pH^i(X) $$
    Then, by part (1) of Theorem~\ref{thm:functorialityofN}, we have $$
  \alpha = \frac{1}{d}f_*f^*(\alpha) \in f_*\left(N^pH^i(X)\right)
  \subset N^pH^i(Y) $$
  since $\dim X = \dim Y$. Hence $f^*(\alpha) \in
  f^* \left(N^pH^i(Y) \right)$.
\end{proof}

The next result is surely known, but we give a proof for lack
of a reference.

\begin{lemma}\label{lemma:basechange}
  Consider a fibered square of varieties
$$
  \xymatrix{
    X'\ar[d]_{g}  \ar[rr]^v && X \ar[d]^f \\
    Y' \ar[rr]_u && Y } $$
  where $u:Y' \to Y$ is an embedding. Then,
  for any $k$ $$
  u^* R^k f_! \Q \simeq R^k g_!( \Q), $$
  where $R^kf_!,\, R^k g_! $ are direct images with proper support.
\end{lemma}
\begin{proof}
  Consider the diagram: $$
  
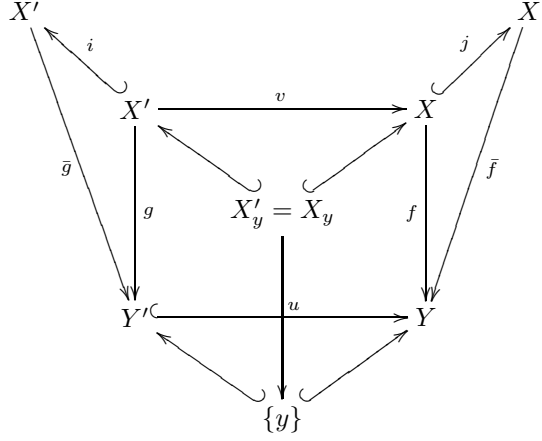
\xymatrix{
    \bar{X'} \ar[dddr]_{\bar{g}} & && & \bar{X} \ar[dddl]^{\bar{f}} \\
    & X' \ar[dd]^g \ar@{_(->}[ul]_i \ar[rr]^v && X \ar[dd]_f \ar@{^(->}[ur]^j  & \\
    && X'_y =X_y \ar@{_(->}[ul] \ar@{^(->}[ur] \ar[dd] && \\
    & Y' \ar@{^(->}[rr]^{\quad u}  && Y &\\
    && \{y\} \ar@{_(->}[ul] \ar@{^(->}[ur] } $$
  where
  $\bar{X}$ (resp. $\bar{X'}$) is a compactification of $X$ (resp. $X'$),
  $i,j$ are inclusions, and $\bar{f}, \bar{g}$ are extensions of $f,g$ to
  compactifications respectively. $X_y$ is the fibre of $f$ 
   at a point $y$ in $Y'$. First, note that we
  have a natural map of sheaves on $Y'$ $$
  u^* R^k \bar{f}_* (j_! \Q)
  \longrightarrow R^k \bar{g}_* (i_! \Q) $$
  Since $\bar{f},\bar{g}$
  are proper morphisms, by \cite[Theorem 6.2]{Iversen}, we have $$
  \left( R^k \bar{f}_* (j_! \Q) \right)_y \cong H^k(\bar{f}^{-1}(y),
  j_! \Q) = H^k_c (f^{-1}(y), \Q) \cong (R^k f_! \Q)_y $$
  and $$
  \left(R^k \bar{g}_* (i_! \Q) \right)_y \cong H^k(\bar{g}^{-1}(y),
  i_! \Q)=H^k_c(g^{-1}(y), \Q) \cong (R^k g_! \Q)_y $$
  Hence we have
  $$
  \left(u^* R^k f_! \Q \right)_y=(R^k f_! \Q)_y = H^k_c
  (g^{-1}(y),\Q) \cong H^k_c(X_y,\Q)=H^k_c(f^{-1},\Q)=(R^k g_! \Q)_y
  $$
  for any $y \in Y'$ and this gives $$
  u^* R^k f_! \Q \simeq R^k
  g_!( \Q) $$

\end{proof}

It will follow, a posteriori, from the theorem that
 the next result holds for all $H$.

\begin{lemma}\label{lemma:hyperplanesection}
  Let $X$ be a smooth projective variety and $X \subset \PP^N$ be a
  fixed embedding. There exists a countable union $B$ of proper
  Zariski closed sets in the dual $(\PP^N)^*$, such that for 
 any hyperplane  $H \notin B$,
  $H^i(X) \to H^i(X \cap H)$ preserves coniveau.
\end{lemma}
\begin{proof}
  Let $\text{Hilb}(X)$ be the Hilbert scheme of $X$ \cite{groth1}, $\mU$ the
  universal family, and $\T=\left(X \times \text{Hilb}(X)\right)-\mU$
  the compliment of $\mU$ in $X \times \text{Hilb}(X)$. Let $p:\T \to
  \text{Hilb}(X)$ the projection. The Hilbert scheme decomposes into
 a countable union $\text{Hilb}(X)=\bigcup^{\infty}_{i=1}C_i$ 
 of irreducible components. Consider the Cartesian diagram
$$
  \xymatrix{
    \T_i \ar[d]_{p_i} \ar[rr] && \T \ar[d]^p \\
    C_i \ar[rr]^{\iota_i} && \text{Hilb}(X) } $$
where  $\iota_i : C_i \to \text{Hilb}(X)$ is the inclusion.
 The fibre $p^{-1}_i(c)\cong X-S_c $, where $S_c$ is a closed subscheme
  of $X$ corresponding to $c \in C_i$.  
  Let $\mR^k=R^kp_! \Q$  and let $\mR^k_i=\iota^*_i \mR^k $ be
  the pull back of $\mR^k$ to $C_i$ for each $i$. Then by
  Lemma~\ref{lemma:basechange}, we have $\mR^k_i=\iota^*_i \mR^k
  \simeq R^k{p_i}_!\Q $.

  First we show that $\mR^k_i$ is a constructible sheaf on $C_i$ for
  each $i$, i.e. there exists a decreasing sequence of Zariski closed
  sets $C_i=F^i_0 \supset F^i_1 \supset \cdots \supset F^i_{n_i}$ such
  that $\mR^k_i|_{F^i_j-F^i_{j+1}}$ is a locally constant sheaf for
  each $j$.  Consider the commutative diagram: $$
  \xymatrix{ \T_i
    \ar[drr]_{p_i} \ar@{^(->}[rr]^{\iota_1} && X \times C_i
    \ar[d]^{\pi}&&(X
    \times C_i)-\T_i=\mU_i \ar[dll]^{h_i} \ar@{_(->}[ll]_{\iota_2} \\
    &&C_i&& } $$
  where $\mU_i$ is the pull back of the universal
  family $\mU$ to $C_i$, $\pi$ is a projection and $\iota_1, \iota_2$
  are inclusions. Note that $\pi$ and $h_i$ are   proper
$\pi$ gives a fibrewise compactification of $p_i$. Hence we
  have an exact sequence of sheaves $$
  \cdots \to R^{k-1}\pi_*\Q
  \stackrel{\rho_{k-1}}{\longrightarrow} R^{k-1}(h_i)_*\Q
  \stackrel{\delta}{\longrightarrow} R^k(p_i)_!\Q
  \stackrel{\gamma_k}{\longrightarrow} R^k \pi_* \Q
  \stackrel{\rho_k}{\longrightarrow} R^k (h_i)_*\Q \to \cdots $$
  From
  this we obtain a short exact sequence $$
  \xymatrix{ 0 \ar[r] &
    \text{Coker}(\rho_{k-1}) \ar[r] & \mR^k_i \ar[r] &
    \text{Ker}(\rho_k) \ar[r] & 0 } $$
  Since for any $k$, $R^k \pi_*
  \Q$ and $R^k(h_i)_* \Q$ are constructible by \cite[Theorem
  2.3.1]{Verdier}, $ \text{Coker}(\rho_{k-1})$ and $\text{Ker}(\rho_k) $
  are also constructible, and consequently so is $\mR^k_i$.  

  Now let $V_{ij}=F^i_j-F^i_{j+1}$
  for each $i,j$. Let 
  $\mW=\{W_i~|~i=1,2,...\}$ be the set of irreducible components of $V_{ij}$.
 Then
  \begin{enumerate}
 \item Each $W_i$ is irreducible.
  \item Each $w \in W_i$ corresponds to a closed subscheme $S_w$ of
    $X$.
  \item $\mR^k|_{W_i}$ is a locally constant sheaf for each $i$
  \end{enumerate}
  Set
  \begin{eqnarray*}
    B_i &=& \{H \in (\PP^N)^*~|~H \supset \bigcup_{w \in W_i}S_w \}\\
    B &=& \bigcup^{\infty}_{i=1} B_i
  \end{eqnarray*}
  Then $B$ is a countable union of proper Zariski closed sets in
  $(\PP^N)^*$. 
 
 Choose $H \in (\PP^N)^* -B$.  We show that
  $H^i (X) \to H^i(X \cap H)$ preserves coniveau. Let
    $\alpha \in N^pH^i(X)$.
   Then there is a subscheme $S$ of $X$ with $\text{codim}(S,X)=q \geq
  p$ such that $$
  \alpha \in \text{Ker}[\psi:H^i(X) \to H^i(X-S)] $$
  Let $[s]$ be the element in $\text{Hilb}(X)$ corresponding to $S$.
  Then $[s] \in W_j$ for some $j$. Since $H \notin B$, $H \notin B_j$,
  equivalently, there is $w \in W_j$ such that $S_{w} \nsubseteq H$.
Set  $k=2 \dim X -i$.
  Since $\mR^k$ is a locally constant sheaf on the connected set $W_j$
  we have an isomorphism $$
  H^k_c (X-S_w ,\Q)
  =\mR^k_w \cong \mR^k_{[s]} = H^k_c (X-S,\Q) $$
  given by parallel transport. Then by duality, we have
  \begin{eqnarray}
    H^i (X-S_w ,\Q) \cong H^i (X-S ,\Q)
  \end{eqnarray}
  Note that $S_w$ is a closed subscheme in  general position with
  respect to $H$, and hence $S_w \cap H$ has codimension 1 in $S_w$.
  Therefore, we have a fibered square $$
  \xymatrix{
    S_w \cap H \ar[d]_{\tau_1} \ar[rr] && X \cap H \ar[d]^{\tau} \\
    S_w \ar[rr] && X } $$
  (where all maps are inclusions) which
  induces a commutative diagram $$
  \xymatrix{
    H^i(X \cap H) \ar[rr]^{\phi \qquad} && H^i((X \cap H) - (S_w \cap H)) \\
    H^i(X) \ar[u]^{\tau^*} \ar[rr]_{\psi_w} && H^i(X-S_w)
    \ar[u]_{\tau^*_1} } $$
  where $\psi_w$ is the composition of $\psi$
  and the isomorphism (1).  Now by the commutativity of above diagram,
  we get $$
  \phi \circ \tau^*(\alpha)= \tau^*_1 \circ \psi_w(\alpha)=0
  $$
  Since $\text{codim}(S_w \cap H, X \cap H)=q \geq p$,
  $$\tau^*(\alpha) \in \text{Ker}[\phi:H^i(X \cap H) \to H^i((X \cap
  H)-(S_w \cap H))]\subset N^pH^i(X \cap H) $$
  Thus, $H^i(X) \to H^i(X
  \cap H)$ preserves coniveau.
\end{proof}

\begin{cor}
  Let $X$ be a smooth projective variety. There exist general hyperplanes
$H_1, H_2,  \cdots , H_p$ such that for $T=X \cap H_1 \cap \cdots \cap H_p$,
  $H^i(X) \to H^i(T)$ preserves coniveau for any $i$.
\end{cor}

\begin{lemma}\label{lemma:commutativity}
  Let $$
  \xymatrix{
    T \ar@{^(->}[rr]^f \ar@{^(->}[d]_{h} &&   S \ar@{^(->}[d]^h \\
    X \ar@{^(->}[rr]_f && Y } $$
  be a fibered square of smooth
  projective varieties such that $f,h$ are embeddings and
  $\text{codim}(T,X)=\text{codim}(S,Y)=c$. Then the following diagram
  commutes: $$
  \xymatrix{
    H^{i-2c}(T) \ar[d]_{h_*} && H^{i-2c}(S) \ar[ll]_{f^*} \ar[d]^{h_*} \\
    H^i(X) && H^i(Y) \ar[ll]^{f^*} } $$
\end{lemma}
\begin{proof}
  We use an argument very similar to the proof of Lemma~\ref{lemma:surjective}. Let
  $N_S$ be a tubular neighborhood of $S$ in $Y$. Let $N_T=N_S \times_Y
  X \cong X \cap N_S$. Then $T \subset N_T$ and $N_T$ is a tubular
  neighborhood of $T$ in $X$, since
  $\text{codim}(T,X)=\text{codim}(S,Y)=c$. Then the  Thom class
$\tau_S$ of $N_S$ restricts to the Thom class $\tau_T$ of $N_T$.
The Gysin map for $S$ (respectively $T$) is given by a composition of 
pull back to $N_S$
(respectively $N_T$), cup product with $\tau_S$ (respectively $\tau_T$) followed
by extension by $0$. Commutativity  of the diagram follows from these remarks.
\end{proof}

\section{End of proof}

Now we give a proof of the third part of the
Theorem~\ref{thm:functorialityofN}

\begin{proof}
  Let $\alpha \in N^p H^i(Y)$ be irreducible, then we have to show that 
$f^*\alpha\in N^pH^i(X)$. By descending induction, we can assume
that $\alpha\notin N^{p+1}H^i(X)$. Therefore $\alpha$ is supported on an
irreducible  subvariety $S\subset Y$ of codimension exactly $p$. Fix a 
desingularization $\tS$ of $S$. Then $\alpha=j_*(\beta)$ for 
some $\beta\in H^{i-2p}(\tS)$, where $j$ is the composition
  of $\tS \to S \hookrightarrow Y$. 

If $f$ is flat and surjective, then we are done by
  Lemma~\ref{lemma:surjective} (the special case suffices). 
Since any $f$ can be factored as
  the inclusion of its graph  followed by a
  projection $X \to X \times Y \to Y$, 
we can reduce to the case where $f$ is a closed
  immersion.

  For the remainder of the proof, assume that $f$ is  an inclusion. 
  Let $d$ be the codimension of $X$ in $Y$. We consider
  three cases.  \newline \newline 
 {\bf Case 1:}  $X$ is not contained in  $S$.\newline

  Note that we may assume that $\text{codim}(S,Y)=p >1$. To see this,
  suppose that $p=1$ and let $Z=X \times_Y S$. Then
  $\text{codim}(Z,X)=1$, since otherwise $X \subset S$. Since $f^*\alpha$
is supported on $Z$, we are done.

  Let $\pi_1:{Y'}=\text{Bl}_{S}Y \to Y$ be the blow up of $Y$
  along $S$, with  ${E'}$ the exceptional divisor. Choose a
desingularization $\pi_2:Y_1 \to {Y'}$ and let $\pi=\pi_1 \circ \pi_2$.
We let $E$ be the strict
  transform of ${E'}$, $X_1$ be the strict transform of $X$ in
  $Y_1$, and $E_1=X_1 \times_{Y_1} E$.
    Then $\text{codim}(E_1,X_1)=1$ by our assumption. First, note that
    by Lemma~\ref{lemma:genericallyfinite}, we have $
    \pi^*(\alpha)$ is supported on $S_1$.
    So, we can find $\beta_1\in H^{i-2p}(\tS_1)$ such that 
${j_1}_* (\beta) =\pi^*(\alpha)$,
    where $j_1$ is the compostion $\tS_1 \to S_1 \to Y_1$ of a
    desingularization of $S_1$ and an inclusion. Now let $H$ be a
    very ample divisor in $Y_1$ and set $$
    S_1 = E \cap H_1 \cap
    \cdots \cap H_{p-1} $$
    where $H_i \in |H|$ are in general position. Note that $\pi(S_1)=S$.
   Let $T=S_1 \times_E E_1 =
    S_1 \cap E_1 = (E\cap H_1 \cap \cdots \cap H_{p-1})\cap E_1=E_1
    \cap H_1 \cap \cdots \cap H_{p-1}$. Hence $T$ is an iterated
    hyperplane section of $E_1$ of codimension $p-1$.  We have 
a commutative diagram
$$
    \xymatrix{
      T=S_1 \times_E E_1 \ar[rr] \ar[d] && S_1 \ar[d]\\
       X_1 \ar[rr]^{f_1} \ar[d]_{\pi} && Y_1 \ar[d]^{\pi}\\
      X \ar[rr]^f && Y } $$
where the top square is Cartesian.
We have $\text{codim}(T,X_1)=\text{codim}(S_1,Y_1)=p$.  
Therefore, by Lemma~\ref{lemma:commutativity} we obtain the
    commutative diagram $$
    \xymatrix{
      H^{i-2p}(T) \ar[d]_{g_*} && H^{i-2p}(S_1) \ar[ll]_{h^*} \ar[d]^{{j_1}_*}\\
      H^i(X_1) && H^i(Y_1) \ar[ll]^{f^*_1} } $$
    So we have $$
    \pi^*
    f^*(\alpha)=f^*_1 \pi^*(\alpha)=f^*_1\circ {j_1}_*(\beta)=g_*
    (h^*(\beta)) \in \text{Im}[H^{i-2p}(T) \to H^i(X_1)] $$
    Hence
    $\pi^* f^*(\alpha) \in N^pH^i(X_1)$ and by
    Lemma~\ref{lemma:genericallyfinite}, $f^*(\alpha) \in N^pH^i(X)$.
    \newline \newline 
   {\bf Case 2:}  $X \subset S$ and $X \neq S$.\\
    
     Then
    $d=\text{codim}(X,Y)>\text{codim}(S,Y)\geq 1$

    Let $\pi: Y' = \text{Bl}_X Y \to Y$ be the blow up of $Y$ along
    $X$, $E$ the exceptional divisor and $S'$ the strict transform of
    $S$. Choose an embedded resolution $Y_1\to Y'$ of singularities,
so that the strict transform $S_1$ of $S$ is nonsingular. Let $E$ be
the strict transform of $E'$.
Let $H$ be a very ample divisor on $Y_1$ and set $$
    X_1 = E    \cap H_1 \cap \cdots \cap H_{d-1} $$
where $H_i\in |H|$ are in general position.
    Note that $\pi(X_1 )=X$.     Let $\iota_1:X_1 \to E$ and
    $\iota_2:E \to Y_1$ be inclusions. Then we have a commutative
    diagram $$
    \xymatrix{ X_1 \ar@{^(->}[rr]^{\iota_1} \ar[drr]_{\pi}
      && E
      \ar@{^(->}[rr]^{\iota_2} \ar[d]^{\pi} && Y_1 \ar[d]^{\pi}\\
      && X \ar@{^(->}[rr]^f && Y } $$
    Let $f_1=\iota_2 \circ \iota_1 :
    X_1 \to Y_1$ be the inclusion. Then the following diagram
    commutes.  $$
    \xymatrix{
      H^i(Y) \ar[rr]^{f^*} \ar@{^(->}[d]_{\pi^*}&& H^i(X) \ar@{^(->}[d]^{\pi^*}\\
      H^i (Y_1 ) \ar[dr]_{\iota^*_2} \ar[rr]^{f^*_1} && H^i(X_1)\\
      &H^i(E) \ar[ur]_{\iota^*_1}& } $$
Recall that we have to show that $f^*\alpha\in N^pH^i(X)$.
By Lemma~\ref{lemma:genericallyfinite}, it suffices to  show that
    $\pi^*(f^*\alpha) \in N^pH^i(X_1)$.
We now chase this element around the other side of the previous
diagram, $$
    \pi^*
    (f^*\alpha)=f^*_1(\pi^*\alpha)=\iota^*_1 \circ
    \iota^*_2(\pi^*\alpha) $$
 Since $\pi:Y_1\to Y$ is a generically
    finite map,  Lemma~\ref{lemma:genericallyfinite} implies that
 $\pi^*(\alpha) \in N^pH^i(Y_1)$. By our assumptions,
 $S_1 \times_{Y_1} E$ has codimension $p$ in $E$, and
$\iota_2^*(\pi^*\alpha)$ is supported on it by 
lemma~\ref{lemma:commutativity}. Therefore $\iota_2^*(\pi^*\alpha)\in N^pH^i(E)$.
  Since $X_1$ is an iterated  hyperplane
    section of $E$ by general hyperplanes, 
    $\iota^*_1$ preserves the coniveau filtration by
    Lemma~\ref{lemma:hyperplanesection}. Hence, 
$\iota^*_1 \circ \iota^*_2(\pi^*\alpha)\in N^pH^i(X_1)$ as required.
 \newline \newline 
 
   {\bf Case 3:} $X =S$.\\

    Then
    $\text{codim}~(X,Y)=\text{codim}~(S,Y)=p$.

    In this case, we have $\beta \in H^{i-2p}(X)$ such that $f_*
    (\beta) =\alpha \in H^i(Y)$. We want to show that $f^*(\alpha) \in
    N^pH^i(X)$. Note that $$
    f^*(\alpha)=f^*f_*(\alpha)=[X]|_{X} \cup
    \alpha $$
    By intersection theory \cite[Section 6.1]{fulton},
    $[X]|_{X}$ is represented by the codimension $p$ algebraic cycle
    $\gamma = [X \cdot X] \in N^pH^{2p}(X)$. So $\gamma =\sum_{k}n_k
    [V_k]$ where $V_k$ is an irreducible subvariety of $X$ of
    codimension $p$ for each $k$. Since $V_k \times X$ is a subvariety
    of $X \times X$ of codimension $p$, we have a Gysin map $$
    \xymatrix{ H^{i-2p}(\tilde{V}_k \times X) \ar[rr]^{\qquad
        g_*}&& H^i(X \times X) } $$
    which sends $\beta \cong 1_k
    \otimes \beta$ to $[V_k] \otimes \beta$, where $1_k$ is a
    generator of $H^0(V_k)$ and $\tilde{V}_k \to V_k$ is a
    desingularization. Also note that $$
    [V_k] \otimes \beta \in
    N^pH^{2p}(X) \otimes N^0H^{i-2p}(X) \subset N^pH^{i}(X \times X)
    $$
    by part (2) of Theorem~\ref{thm:functorialityofN}. Let $\Delta:X \to
    X \times X$ be the diagonal map. Then, by letting $S_k=V_k \times
    X$, $Y_k=X \times X$, and $X_k=X$, for each $k$, we can reduce
    this Case to the Case 1. i.e. we have $$
    \Delta^* : H^i(Y_k)
    \longrightarrow H^i(X_k) $$
    with $$
    [V_k] \otimes \beta \in
    \text{Im}[H^{i-2p}(\tilde{S}_k) \to H^i(Y_k)] \subset
    N^pH^i(Y_k) $$
    and $S_k \neq X_k$, where $\tS_k \to S_k$ is a
    desingularization of $S_k$. Then we have $$
    \Delta^*([V_k] \otimes
    \beta)=[V_k] \cup \beta \in N^pH^i(X_k)=N^pH^i(X) $$
    for each $k$.
    Therefore, $$
    f^*(\alpha)=\sum_k n_k \left([V_k] \cup \beta
    \right) \in N^pH^i(X) $$
    This completes the proof of the Theorem.

  \end{proof}

\end{document}